\documentclass[11pt]{article}
\usepackage{amsfonts}
\usepackage{mathrsfs}
\usepackage{amsmath}
\usepackage{amssymb}
\usepackage[mathscr]{eucal}
\usepackage{graphicx}
\renewcommand{\paragraph}{\roman{paragraph}}

\def \ov{\overline}
\def \C{\mathcal{C}}
\def \e{\epsilon}
\newtheorem{theorem}{\scshape \mdseries  Theorem}[section]
\newtheorem{lemma}[theorem]{\scshape \mdseries  Lemma}
\newtheorem{coro}[theorem]{\scshape \mdseries  Corollary}

\begin{document}

\title{\sf Spectral conditions for a graph to be Hamilton-connected }
\author{Gui-Dong Yu$^{1,2,}$\thanks{Email: yuguid@aqtc.edu.cn.
Supported by NSF of Department of Education of Anhui Province (KJ2011A195)
and Innovation Fund  for Graduates of Anhui University.},\ \
Yi-Zheng Fan$^{1,}$\thanks{Corresponding author.
Email: fanyz@ahu.edu.cn.
Supported by National Natural Science Foundation of China (11071002),
Program for New Century Excellent Talents in University, Key Project of Chinese Ministry of Education (210091),
Specialized Research Fund for the Doctoral Program of Higher Education (20103401110002),
Science and Technological Fund of Anhui Province for Outstanding Youth  (10040606Y33)
Project of Educational Department of Anhui Province (KJ2010B136),
Scientific Research Fund for Fostering Distinguished Young Scholars of Anhui University (KJJQ1001),
Project for Academic Innovation Team of Anhui University (KJTD001B).} \\
  {\small  \it $1.$ School of Mathematical Sciences, Anhui University, Hefei 230039, P.R. China}\\
    {\small  \it $2.$ School of Mathematics \& Computation Sciences, Anqing Normal College, Anqing 246011, P.R. China}
  }
\date{}
\maketitle

\noindent {\bf Abstract:}
In this paper we establish some
spectral conditions for a graph to be Hamilton-connected in terms of
the spectral radius of the adjacency matrix or the signless Laplacian of
the graph or its complement. For the existence of Hamiltonian paths or cycles in a graph,
we also give a sufficient condition by the signless Laplacian spectral radius.


\noindent {\bf MR Subject Classifications:}  05C50, 05C45, 05C35

\section{Introduction}
Let $G=(V,E)$ be a simple graph of order $n$
   with vertex set $V=V(G)=\{v_1,v_2,\ldots,v_n\}$ and edge set $E=E(G)$.
The {\it adjacency matrix} of $G$ is defined to be a matrix $A(G)=[a_{ij}]$ of order $n$,
  where $a_{ij}=1$ if $v_{i}$ is adjacent to $v_{j}$, and $a_{ij}=0$ otherwise.
  The largest eigenvalue of $A(G)$, denoted by $\mu(G)$, is the spectral radius of $A(G)$.
The degree matrix of $G$ is denoted by $D(G)=\hbox{diag}\left(d_G(v_1),d_G(v_2), \ldots, d_G(v_n)\right)$,
   where $d_G(v)$ or simply $d(v)$ denotes  the degree of a vertex $v$ in the graph $G$.
The matrix $L(G)=D(G)-A(G)$ is the {\it Laplacian matrix} of $G$, and
the matrix $Q(G)=D(G)+A(G)$ is the {\it signless Laplacian matrix} (or {\it $Q$-matrix}) of $G$.
Denote by $\gamma(G)$ the largest eigenvalues (or the  spectral radius) of $Q(G)$.

A  {\it Hamiltonian cycle} of the graph $G$ is a cycle of order $n$
contained in $G$, and a {\it Hamiltonian path} of $G$ is a path of
order $n$ contained in $G$. A graph $G$ is said to be {\it
Hamiltonian} if it contains Hamiltonian cycles. The problem of
deciding whether a graph is Hamiltonian is one of
   the most difficult classical problems in graph theory.
Indeed, determining whether a graph is Hamiltonian is NP-complete.

Recently, the spectral theory of graphs has been applied to this
problem. Fiedler and Nikiforov \cite{fied}  give sufficient
conditions for a graph having Hamiltonian paths or cycles
  in terms of the spectral radius of the adjacency matrix of the graph or its complement.
Zhou \cite{zhou} investigates the spectral radius of the signless Laplacian matrix of the complement of a graph,
 and present some conditions for the existence of Hamiltonian cycles or paths.
Butler and Chung \cite{but} establish a sufficient condition for a graph being Hamiltonian, that is,
the nontrivial eigenvalues of the Laplacian are sufficiently close to the average degree of the graph.
Their result does not need the graph to be dense, but it implies the graph is almost regular.

A graph $G$  is  called {\it Hamilton-connected}   if every two
vertices of $G$ are connected by a Hamiltonian path. Surely all
Hamilton-connected graphs are Hamiltonian.
In this paper, motivated by the ideas in \cite{fied} and \cite{zhou}, we
establish the spectral conditions for a graph to be
Hamilton-connected in terms of the spectral radius of the adjacency
matrix or signless Laplacian of the graph or its complement. In
addition, we also give the condition on the signless Laplacian
spectral radius of a graph for the existence of Hamiltonian paths or
cycles.

\section{Preliminaries}
Let $K_{n}$ be the complete graph on $n$ vertices.
Write $K_{n-1}+v$ for $K_{n-1}$ together with an isolated vertex $v$,
   $K_{n-1}+e$ for $K_{n-1}$ together with a pendent edge $e$,
   and $K_{n-1}+e+e'$ for $K_{n-1}$ together with a vertex joining two vertices of $K_{n-1}$
   by the edges $e,e'$, respectively.
Denote by $K_{n,m}$  a complete bipartite graph with two parts
having $n,m$ vertices, respectively. The graph $K_{1,n-1}$ is also
called a star on $n$ vertices.
Denote by $\e(G)$ the number of edges of a graph $G$.

Ore \cite{ore} showed that
if  $$d(u)+d(v)\geq n-1, \eqno(2.1)$$
for any pair of nonadjacent vertices $u$ and $v$,
then $G$ contains a Hamiltonian path.
If the inequality (2.1) is strict, then $G$ contains a Hamiltonian cycle.
Erdo\"s and Gallai \cite{erdos} proved that if
$$d(u)+d(v)\geq n+1, \eqno(2.2) $$
for any pair of nonadjacent vertices $u$ and $v$, then $G$ is
Hamilton-connected. So we have the following facts, where the first two facts are due to
Fiedler and Nikiforov \cite{fied}.

\begin{lemma} \label{fied}
Let $G$ be a graph of order $n$ with $m$ edges.
Assume that $$m = \binom{n-1}{2} +r. $$
If $r \ge 0$, then $G$ contains a Hamiltonian path unless $G=K_{n-1}+v$.
If $r \ge 1$, then $G$ contains a Hamiltonian cycle unless $G=K_{n-1}+e$.
If $r \ge 2$, then $G$ is Hamilton-connected unless $G=K_{n-1}+e+e'$.
\end{lemma}

For an integer $k \geq 0$, the {\em $k$-closure of a graph} $G$, denoted by $\C_k(G)$,
   is the graph obtained from $G$ by successively joining pairs of nonadjacent vertices
   whose degree sum is at least $k$ until no such pair remains; see \cite{bondy}.
The $k$-closure of the graph $G$ is unique, independent of the order in which edges are added.
Note that
$d_{\C_{k}(G)}(u)+d_{\C_{k}(G)}(v)\leq k-1$ for any pair of nonadjacent vertices $u$ and $v$ of $C_{k}(G)$.

\begin{lemma} {\em \cite{bondy}} \label{clo}
{\em(i)}  A graph $G$ has a Hamilton path if and only if $\C_{n-1}(G)$ has one.

{\em(ii)}  A graph $G$ has a Hamilton cycle if and only if $\C_{n}(G)$ has one.

{\em(iii)} A graph $G$ is Hamilton-connected if and only if $\C_{n+1}(G)$ is so.

\end{lemma}

\begin{lemma} {\em \cite{tom}} \label{tom}
Any $t$-regular graph of order $2t \;(t \geq3)$ not isomorphic to $K_{t,t}$,
or of order $2t+1$ for even $t \geq 4)$,  is Hamilton-connected.
\end{lemma}

At the final of this section, we introduce some bounds for the
spectral radius of the adjacency matrix or the signless Laplacian
matrix of a graph. Let $G$ be a graph containing a vertex $v$.
Denote $m_G(v)=m(v)={1 \over d_G(v)} \sum_{u \in N_G(v)} d_G(u)$ if $d_G(v) >0$, and
$m_G(v)=0$ otherwise, where $N_G(v)$ or simply $N(v)$ denotes the neighborhood of $v$ in $G$.

\begin{lemma} {\em \cite{sta}} \label{sta}
Let $G$ be graph with $m$ edges.
Then
$$ \mu(G) \le -\frac{1}{2}+\sqrt{2m+\frac{1}{4}},$$
with equality if and only if $G$ is a complete graph together with some isolated vertices.
\end{lemma}

\begin{lemma} {\em \cite{das1}} \label{avedg}
Let $G$ be a graph of order $n$ with $m$ edges.
Then
$$ \max\{d(v)+m(v): v \in V(G)\} \le \frac{2m}{n-1}+n-2,$$
with equality if and only if $G \supseteq K_{1,n-1}$ or $G=K_{n-1}+v$.
\end{lemma}

\begin{lemma}{\em \cite{das2}} \label{upbd}
Let G be a connected graph. Then,
$$ \gamma(G) \le  \max\{d(v)+m(v): v \in V(G)\},\eqno(2.3)$$
with equality if and only if $G$ is a regular graph or a bipartite semiregular graph.
\end{lemma}

In Lemma \ref{upbd}, if $G$ is disconnected, there exists a  component $G'$ of $G$ such that
$$\gamma(G)=\gamma(G')\le  \max\{d_{G'}(v)+m_{G'}(v): v \in V(G')\}
\le  \max\{d_G(v)+m_G(v): v \in V(G)\}.$$
So the inequality (2.3) also holds.
Combining Lemma \ref{avedg} and Lemma \ref{upbd}, we have the following result; also see \cite{feng}.

\begin{coro} \label{upbd2}
Let $G$ be a graph of order $n$ with $m$ edges.
Then
$$\gamma(G) \le \frac{2m}{n-1}+n-2. \eqno(2.4)$$
If $G$ is connected, then the equality in (2.4) holds if and only if $G=K_{1,n-1}$ or $G=K_n$.
Otherwise, the equality in (2.4) holds if and only if $G=K_{n-1}+v$.
\end{coro}

\section{Main Results}

\begin{theorem}
Let $G$ be a graph of order $n$. If
$$\mu(G) > -\frac{1}{2}+\sqrt{\left(n-\frac{3}{2}\right)^2+2}, \eqno(3.1)$$
then $G$ is Hamilton-connected unless $G=K_{n-1}+e+e'$.
\end{theorem}

{\bf Proof:}
Suppose $G$ has $m$ edges.
By Lemma \ref{sta} together with (3.1),
$$ -\frac{1}{2}+ \sqrt{\left(n-\frac{3}{2}\right)^2+2}< \mu(G) \le -\frac{1}{2}+\sqrt{2m+\frac{1}{4}}.$$
Hence $m \ge \binom{n-1}{2}+2$, and the result follows
by Lemma \ref{fied}. \hfill $\blacksquare$

\begin{theorem}
Let $G$ be a graph of order $n \ge 4$. If
$$\mu(\ov{G}) < \sqrt{\frac{(n-2)^2}{n}}, \eqno(3.2)$$
then $G$ is Hamilton-connected.
\end{theorem}

{\bf Proof:} Let $H=\C_{n+1}(G)$. If $H=K_{n}$, then the result
follows from Lemma \ref{clo}. Suppose that $H \neq K_{n}$ and $G$ is
not Hamilton-connected. Then $H$ is also not Hamilton-connected by
Lemma \ref{clo}. Observe that $d_{H}(u)+d_{H}(v)\leq n$ for any pair
of nonadjacent vertices $u$ and $v$ (always existing) in $H$. Thus,
for any edge $uv \in E(\ov{H})$,
$d_{\ov{H}}(u)+d_{\ov{H}}(v)=2(n-1)-[d_{H}(u)+d_{H}(v)]\geq n-2$. It
follows that
$$\sum_{u\in V(\ov{H})}d_{\ov{H}}^{2}(u)
 =\sum_{uv\in E(\ov{H})}[d_{\ov{H}}(u)+d_{\ov{H}}(v)]\geq (n-2)\e(\ov{H}).$$
Using the inequality of Hofmeister \cite{hof},
$$n \mu(\ov{H})^2 \ge \sum_{u\in V(\ov{H})}d_{\ov{H}}^{2}(u) \ge (n-2)\e(\ov{H}).$$
As $\ov{H} \subseteq \ov{G}$, by Perron-Frobenius theorem and (3.2), we have
$ \mu(\ov{H}) \le \mu(\ov{G}) < \sqrt{\frac{(n-2)^2}{n}}$,
and hence
$$ (n-2)^2 > n \mu(\ov{G})^2 \ge n \mu(\ov{H})^2 \ge (n-2) \e(\ov{H}).$$
So, $\e(\ov{H}) < n-2$ and $\e(H) \ge \binom{n}{2} -(n-3)=\binom{n-1}{2} +2.$ Since $H$ is not
Hamilton-connected, by Lemma \ref{fied}, $H=K_{n-1}+e+e'$. Note that
$\ov{H}=K_{1,n-3} \oplus \ov{K_2}$. So, $\ov{G}$ contains
$K_{1,n-3}$ as a subgraph, and $\mu(\ov{G}) \ge \sqrt{n-3} \ge
\sqrt{\frac{(n-2)^2}{n}}$, a contradiction. \hfill $\blacksquare$

\begin{theorem}
Let $G$ be a graph of order $n$. If
$$\gamma(G) > 2(n-2)+\frac{2}{n-1}, \eqno(3.3)$$
then $G$ is Hamilton-connected  unless $G=K_{n-1}+e+e'$.
\end{theorem}

{\bf Proof:}
Suppose $G$ has $m$ edges.
By Corollary \ref{upbd2} together with (3.3),
$$ 2(n-2)+\frac{2}{n-1} < \gamma(G) \le \frac{2m}{n-1}+n-2.$$
Hence $m \ge \binom{n-1}{2} +2$, and the result follows by Lemma \ref{fied}. \hfill $\blacksquare$

Let $G_1,G_2$ be two vertex-disjoint graphs.
Denote by $G_1 \oplus G_2$ the union of $G_1,G_2$.
Denote by $G_1 \vee G_2$ the join of $G_1$ and $G_2$,
  a graph from $G_1 \oplus G_2$ by adding all possible edges between each vertex of $G_1$ and each vertex of $G_2$.
 Let $\mathbb{ES}_n$ be the set of following graphs of order $n$:
$(K_s \oplus K_{n-2-s})\vee K_{2}$, $K_{\frac{n}{2},\frac{n}{2}}$,  and $R \vee K_{r}$,
 where $n$ is even for the latter two graphs,
 $R$ is a regular graph of order $n-r$ with degree $\frac{n}{2}-r$,
 $1 \le s \le n-3$, $1\leq r \leq \frac{n}{2}$.

\begin{theorem}
Let $G$ be a graph of order  $n \ge 6$, and $G \notin \mathbb{ES}_n$. If
$$\gamma(\ov{G})\leq n-2,  \eqno(3.4)$$
then $G$ is Hamilton-connected.
\end{theorem}

{\bf Proof:} Let $H=\C_{n+1}(G)$. If $H=K_{n}$, then the result
follows from Lemma \ref{clo}. Suppose that $H \neq K_{n}$ and $G$ is
not Hamilton-connected. Then $H$ is also not Hamilton-connected by
Lemma \ref{clo}. Observe that $d_{H}(u)+d_{H}(v)\leq n$ for any pair
of nonadjacent vertices $u$ and $v$ (always existing) in $H$. Thus,
for any edge $uv \in E(\ov{H})$,
$d_{\ov{H}}(u)+d_{\ov{H}}(v)=2(n-1)-[d_{H}(u)+d_{H}(v)]\geq n-2$. It
follows that
$$Z(\ov{H}):=\sum_{u\in V(\ov{H})}d_{\ov{H}}^{2}(u)
 =\sum_{uv\in E(\ov{H})}(d_{\ov{H}}(u)+d_{\ov{H}}(v))\geq (n-2)\e(\ov{H}). \eqno(3.5)$$
By Lemma 3 of \cite{zhou},
$$ \gamma(\ov{H}) \ge \frac{Z(\ov{H})}{\e(\ov{H})}. \eqno(3.6)$$
As $\ov{H}$ is a subgraph of $\ov{G}$, by Perron-Frobenius theorem, $\gamma(\ov{G}) \ge \gamma(\ov{H})$.
Thus, by (3.4), (3.5) and (3.6),
$$ n-2 \ge \gamma(\ov{G}) \ge \gamma(\ov{H}) \ge \frac{Z(\ov{H})}{\e(\ov{H})} \ge n-2,$$
and hence all the inequalities become equalities.
So, $d_{\ov{H}}(u)+d_{\ov{H}}(v)=n-2$ for any edge $uv\in E(\ov{H})$ by  the fact mentioned in \cite{zhou}.
This implies that $\ov{H}=K_{\frac{n}{2}}\cup K_{\frac{n}{2}}$ for an even $n$
or $\ov{H}$ contains exactly one nontrivial component $F$,
where $F$ is either regular or bipartite semi-regular and $\frac{n}{2} \leq|V(F)|\leq n$.
Since $\gamma(\ov{G}) = \gamma(\ov{H})$ and $\ov{H}$ is a subgraph of $\ov{G}$, by Perron-Frobenius theorem
we have $\ov{G}=\ov{H}$ and then $G=H$.
If $\ov{H}=K_{\frac{n}{2}}\cup K_{\frac{n}{2}}$, then
$G=H=K_{\frac{n}{2},\frac{n}{2}}$, which contradicts the assumption on $G$.

Now suppose $\ov{H}$ contains exactly one nontrivial component $F$.
First assume $F$ is bipartite semi-regular with two parts $U_1,U_2$.
Let $u_1u_2$ be an edge of $F$,  where $u_1 \in U_1, u_2 \in U_2$.
Then $V(F)  \supseteq (N_F(u_1) \cup N_F(u_2))$, which implies $F$ contains at least $n-2$ vertices.

If $F$ contains exactly $n-2$ vertices, then $F=K_{s,n-2-s}$ for some positive integer $s$,
  and $\ov{G}=\ov{H}=K_{s,n-2-s} \oplus \ov{K_2}$.
So $G=(K_s \oplus K_{n-2-s})\vee K_{2}$, a contradiction to the assumption on $G$.

If $F$ contains exactly $n-1$ vertices, then $U_1$ or $U_2$, say $U_1$,
  contains an additional vertex rather than those in $N(u_1) \cup N(u_2)$,
then by summing the degrees of the vertices in $U_1$ and $U_2$ respectively, we have
$(d_{\ov{H}}(u_2)+1)d_{\ov{H}}(u_1)=d_{\ov{H}}(u_1) d_{\ov{H}}(u_2)$, a contradiction.

If $F$ contains $n$ vertices, or equivalently $\ov{H}=F$ is connected,
then by a similar discussion as the above,
  $U_1$ (or $U_2$) consists of the vertices of $N(u_2)$ (or $N(u_1)$) and an additional vertex.
As $(d_{\ov{H}}(u_2)+1)d_{\ov{H}}(u_1)=(d_{\ov{H}}(u_1)+1)d_{\ov{H}}(u_2)$,
we have $d_{\ov{H}}(u_1)=d_{\ov{H}}(u_2)=\frac{n-2}{2}$, which implies $\ov{H}$ is a bipartite regular graph.
Hence, $\ov{G}=\ov{H}=K_{\frac{n}{2},\frac{n}{2}}-M$,
  where $M$ is a perfect matching of $K_{\frac{n}{2},\frac{n}{2}}$.
It is easy to check that $G$ is Hamilton-connected.

Finally we assume $F$ is a regular graph.
Then $F$ is regular of degree $\frac{n}{2}-1$ for an even $n$,
and $F$  is not bipartite semi-regular.
If $F=\ov{H}$, then $\ov{G}=\ov{H}$, and $G$ is regular of degree $\frac{n}{2}$.
Because $G \ne K_{\frac{n}{2},\frac{n}{2}}$, by Lemma \ref{tom}, $G$ is Hamilton-connected.
Otherwise,  $\ov{H} \supsetneq F$, and hence $\ov{G}=\ov{H}= F \oplus \ov{K_{r}}$,
   where $r=n-|V(F)|$ and $1\leq r\leq \frac{n}{2}$.
So, $G=\ov{F} \vee K_{r}$, where $\ov{F}$ is regular of order $n-r$ and degree $\frac{n}{2}-r$,
a contradiction to the assumption on $G$. The result follows. $\blacksquare$

\section{Remark}
Zhou \cite{zhou} gives a sufficient condition for a graph having
Hamiltonian paths or cycles in terms of the signless Laplacian
spectral radius of the complement of the graph. Here we give a
complementary condition for the existence of Hamiltonian paths or
cycles in terms of the signless Laplacian spectral radius of the
graph.

\begin{theorem} \label{sgnH}
Let $G$ be a graph of order $n$.
If $$\gamma(G)\geq 2(n-2), \eqno(4.1) $$
then $G$ contains a Hamiltonian path unless $G=K_{n-1}+v$.
If the inequality (4.1) is strict, then $G$ contains a Hamiltonian cycle unless $G=K_{n-1}+e$.
\end{theorem}

{\bf Proof:}
Suppose $G$ has $m$ edges.
By Corollary \ref{upbd2} and the assumption,
we have
$$2(n-2)\leq \gamma(G)\leq \frac{2m}{n-1}+n-2.$$
Hence we obtain $m \geq \binom{n-1}{2} $, with strict inequality if (4.1) is strict.
The result follows by lemma \ref{fied}. \hfill $\blacksquare$

Fiedler and Nikiforov \cite{fied} gave a condition for the existence
of Hamiltonian paths or cycles by the spectral radius of the
adjacency matrix of a graph.

\begin{theorem} {\em \cite{fied} }\label{adjH}
Let $G$ be a graph of order $n$.
If $$\mu(G)\geq n-2,\eqno(4.2)$$
then $G$ contains a Hamiltonian path unless $G=K_{n-1}+v$.
If the inequality (4.2) is strict,
then $G$ contains a Hamiltonian path unless $G=K_{n-1}+e$.
\end{theorem}

Note that $\gamma(G) \ge 2 \mu(G)$ (see \cite{chen}).
So, if $\mu(G) \ge n-2$, then $\gamma(G) \ge 2(n-2)$, and hence
$G$ admits a Hamiltonian path, or Hamiltonian cycle if the equality holds strictly.
However, there exist graphs $G$ with $\gamma(G) \ge 2(n-2)$ but $\mu(G) < n-2$, which contains
a Hamiltonian path or cycle.

For example, let $G=(K_r \oplus K_r) \vee K_s$ of order $n:=2r+s$,
where $r \ge 2$ and $2(r-1)^2 \le s < (2r-1)(r-1)$. Surely, the
graph $G$ contains Hamiltonian cycles. By a little computation,
$\mu(G)$ is the largest root of the polynomial
$f(x)=[x-(r-1)][x-(s-1)]-2rs$ and $\gamma(G)$ is the largest root of
the polynomial $g(x)=[x-(2r+s-2)][x-(2r+2s-2)]-2rs$. Observe that
$f(n-2)=f(2r+s-2)=(2r-1)(r-1)-s>0$, so that $\mu(G)<n-2$, and
$g(2(n-2))=g(2(2r+s-2))=4(r-1)^2-2s \le 0$, so that $\gamma(G) \ge
2(n-2)$. So, we can apply Theorem \ref{sgnH} but not Theorem
\ref{adjH} for the existence of Hamiltonian paths or cycles of $G$.

\end{document}